\documentclass[11pt]{article}
\usepackage[T2A]{fontenc}
\usepackage{physics}
\usepackage{amsmath}
\usepackage{tikz}
\usepackage{mathdots}
\usepackage{yhmath}
\usepackage{cancel}
\usepackage{color}
\usepackage{siunitx}
\usepackage{array}
\usepackage{multirow}
\usepackage{amssymb}
\usepackage{gensymb}
\usepackage{tabularx}
\usepackage{extarrows}
\usepackage{booktabs}
\usetikzlibrary{fadings}
\usetikzlibrary{patterns}
\usetikzlibrary{shadows.blur}
\usetikzlibrary{shapes}
\usepackage[english]{babel}
\usepackage{mathbbol}
\usepackage{amsthm,amssymb}

\usepackage{hyperref}
\usepackage{enumerate}

\title{The zeta function of zeros and poles of a meromorphic function, the criterion for the absence of zeros, and its application to sums of Cauchy kernels}
\author{Vladimir Shemyakov \\
St. Petersburg, Russia \\
	\texttt{vladimir.v.shemyakov@gmail.com}
}
\date{\today}

\begin{document}
\maketitle
\newtheorem{Hyp}{Conjecture}
\newtheorem{thm}{Theorem}
\newtheorem{definition}{Definition}
\newtheorem{remark}{Remark}
\newtheorem{lema}{Lemma}

\begin{abstract}
   A connection between the zeta functions of zeros and poles of a meromorphic function has been established, and using it, a criterion for the absence of zeros has been derived. Sufficient conditions for the existence of zeros of sums of Cauchy kernels have been obtained, including those based on geometric conditions on the poles.
\end{abstract}
\section{Introduction}

Let $\displaystyle f$ be a meromorphic function in $\displaystyle \mathbb{C}$ with poles $\displaystyle \{t_{n}\}_{n\in \mathbb{N}} \subset \mathbb{C}$, and let $\displaystyle G_{n}( z)$ be the principal part of the Laurent series of the function $\displaystyle f$ at the point $\displaystyle t_{k}$:
\begin{equation*}
G_{n}( z) =\sum\limits _{k=1}^{N_{n}}\frac{a_{k}^{n}}{( z-t_{n}){^{k}}}.
\end{equation*}
Then, according to the Mittag-Leffler theorem, there exists a sequence of polynomials $\displaystyle P_{n}$ such that
\begin{equation*}
f=\sum\limits _{n=1}^{\infty }( G_{n} -P_{n}).
\end{equation*}
such that the convergence is uniform on compact sets not containing the points $\displaystyle \{t_{n}\}_{n\in \mathbb{N}}$. In other words, the theorem describes the structure of an arbitrary meromorphic function, albeit not very explicitly. In analysis, a significant subclass of meromorphic functions plays a distinct role, represented in the form of series
\begin{equation*}
\sum\limits _{k=1}^{\infty }\frac{c_{k}}{z-t_{k}}.
\end{equation*}
These functions are called sums of Cauchy kernels, with each term $\displaystyle \frac{c_{k}}{z-t_{k}}$ being a Cauchy kernel. The term "kernel" has a specific meaning; it is indeed the kernel of the integral operator in the Cauchy integral formula, which reconstructs the value of an analytic function inside a region $\displaystyle \Omega $ with a rectifiable boundary from the values on the boundary $\displaystyle \partial \Omega $
\begin{equation*}
g( z) =\oint\limits _{\partial \Omega }\frac{f( \xi )}{\xi -z} d\xi ,\ z\in \Omega .
\end{equation*}
This theorem essentially states that an analytic function is a "continuous" sum of Cauchy kernels, if one thinks of the integral as a sum. Here are some well-known examples of sums of Cauchy kernels:
\begin{equation*}
\sum\limits _{n=-\infty }^{\infty }\frac{(-1)^{n}}{z-n} =\frac{\pi }{\sin( \pi z)} ,\ \ \sum _{n=-\infty }^{\infty }\frac{( -1)^{n}}{n(z-n)} =\frac{\pi }{z\sin\left( \pi z\right)}-\frac{1}{z^2}  \ \ -\ \text{no zeros},
\end{equation*}
\begin{equation*}
\sum\limits _{n=1}^{\infty }\frac{1}{k(z-n)} =\frac{\gamma +\Psi (1-z)}{z} ,\ \ \sum\limits _{n=1}^{\infty }\frac{1}{n^{2}(z-n)} =\frac{\pi ^{2}}{6} +\frac{\gamma +\Psi ( 1-z)}{z^{2}} \ -\  \text{infinitely many zeros},
\end{equation*}
\begin{equation*}
\sum\limits _{n=1}^{\infty }\frac{1}{z-n^{2}} =\frac{\pi \cot\left( \pi \sqrt{z}\right)\sqrt{z} -1}{2z} ,\ \ \ \sum\limits _{n=1}^{\infty }\frac{1}{2^{n}( z-n)} =-\frac{\Phi \left(\frac{1}{2} ,1,-z\right) z+1}{z} \ -\ \text{infinitely many zeros}, 
\end{equation*}
where $\displaystyle \Psi =\frac{\Gamma '}{\Gamma }$, $\displaystyle \Gamma $ is the Euler's Gamma function, and $\displaystyle \Phi $ is the Lerch's zeta function.
For further convenience, let's present a classical result guaranteeing meromorphy in $\displaystyle \mathbb{C}$ functions of the form $\displaystyle \sum _{k=1}^{\infty }\frac{c_{k}}{( z-t_{k})^{m}}$.
\begin{lema}\label{lem:classic}
Let $\displaystyle \{c_{n}\}_{n\in \mathbb{N}} ,\{t_{n}\}_{n\in \mathbb{N}} \subset \mathbb{C}$, such as $\displaystyle \lim _{n\rightarrow \infty } |t_{n} |=\infty $ and
\begin{equation*}
\sum _{n=1}^{\infty }\frac{|c_{n} |}{|t_{n} |^{m}} < +\infty \ ,\ m\in \mathbb{N} .
\end{equation*}
Then the function $\displaystyle \sum _{k=1}^{\infty }\frac{c_{k}}{( z-t_{k})^{m}}$ is meromorphic in $\displaystyle \mathbb{C}$, the series converges absolutely, uniformly on compacts in $\displaystyle \mathbb{C}$ not containing points $\displaystyle \{t_{n}\}_{n\in \mathbb{N}}$, and it can be differentiated term by term infinitely many times.
\end{lema}
Therefore, throughout, we will require this condition from Lemma \ref{lem:classic}, but thanks to the example
\begin{equation*}
\sum _{k=-\infty }^{\infty }\frac{( -1)^{k}}{z-k} =\frac{\pi }{\sin( \pi z)}
\end{equation*}
и $\displaystyle \sum _{k=-\infty }^{\infty }\left| \frac{( -1)^{k}}{k}\right| =+\infty $ we see that the condition in Lemma \ref{lem:classic} is not necessary. Looking at the examples above, a logical hypothesis comes to mind
\begin{Hyp}
[Clunie, Eremenko, Rossi,\cite{ClunieEremenkoRossi}, conj. 2.7]\label{hyp:ClunieEremenkoRossi}

Let \textcolor[rgb]{0.82,0.01,0.11}{$\displaystyle \{c_{n}\}_{n\in \mathbb{N}} \geqslant 0$} and $\displaystyle \{t_{n}\}_{n\in \mathbb{N}} \subset \mathbb{C}$ such that $\displaystyle \lim _{n\rightarrow \infty } |t_{n} |=\infty $ and
\begin{equation*}
\sum\limits _{n=1}^{\infty }\frac{c_{n}}{|t_{n} |} < +\infty.
\end{equation*}
Then the meromorphic function in $\displaystyle \mathbb{C}$ given by $\displaystyle \sum\limits _{k=1}^{\infty }\frac{c_{k}}{z-t_{k}}$ has infinitely many zeros in $\displaystyle \mathbb{C}$.
\end{Hyp}
The examples do not contradict the hypothesis; however, it has not been proven or disproven yet. Nonetheless, there are various results in this direction, such as the classical
\begin{thm}
[\cite{ClunieEremenkoRossi}, th. 2.1]

Let \textcolor[rgb]{0.82,0.01,0.11}{$\displaystyle \{c_{n}\}_{n\in \mathbb{N}} \subset \mathbb{N}$} and $\displaystyle \{t_{n}\}_{n\in \mathbb{N}} \subset \mathbb{C}$ such that $\displaystyle \lim _{n\rightarrow \infty } |t_{n} |=\infty $ and
\begin{equation*}
\sum\limits _{n=1}^{\infty }\frac{c_{n}}{|t_{n} |} < +\infty.
\end{equation*}
Then, the meromorphic function in $\displaystyle \mathbb{C}$ given by $\displaystyle \sum\limits _{k=1}^{\infty }\frac{c_{k}}{z-t_{k}}$ has infinitely many zeros in $\displaystyle \mathbb{C}$.
\end{thm}

\begin{thm}
[Clunie, Eremenko, Rossi,\cite{ClunieEremenkoRossi}, th. 2.10]

Let \textcolor[rgb]{0.82,0.01,0.11}{$\displaystyle \{c_{n}\}_{n\in \mathbb{N}} \geqslant 0$} and $\displaystyle \{t_{n}\}_{n\in \mathbb{N}} \subset \mathbb{C}$ such that $\displaystyle \lim _{n\rightarrow \infty } |t_{n} |=\infty $ and
\begin{equation*}
\sum\limits _{n=1}^{\infty }\frac{c_{n}}{|t_{n} |} < +\infty ,
\end{equation*}
\begin{equation*}
\sum\limits _{|t_{n} |< r} c_{n} =o\left(\sqrt{r}\right) \ \text{at} \ r\rightarrow \infty .
\end{equation*}
Then, the meromorphic function in $\displaystyle \mathbb{C}$ given by $\displaystyle \sum\limits _{k=1}^{\infty }\frac{c_{k}}{z-t_{k}}$ has infinitely many zeros in $\displaystyle \mathbb{C}$.
\end{thm}

There are results aimed at removing restrictions from $\displaystyle c_{n}$ while imposing more requirements on $\displaystyle t_{n}$, for example, bounding the growth of $\displaystyle t_{n}$ from above and below. To achieve this, let's introduce standard characteristics from Nevanlinna theory - a powerful branch of complex analysis that uniformly studies the distribution of values of meromorphic functions. Herman Weyl referred to Rolf Nevanlinna's theory as "one of the greatest mathematical events of the twentieth century".

\begin{definition}[R. Nevanlinna, \cite{Book_Gold_Ostr}]
Let $\displaystyle f$ be meromorphic in $\displaystyle \mathbb{C}$ and $\displaystyle \{t_{n}\}_{n\in \mathbb{N}} \subset \mathbb{C}$ be the poles of $\displaystyle f$ with multiplicity taken into account. Let's introduce
\begin{equation*}
n( r,f) =\ \sum\limits _{|t_{k} |< r} 1=\#\{\ k\ |\ |t_{k} |< r\} \ 
\end{equation*}
- the number of poles inside the disk $\displaystyle |z|\leqslant r$,
\begin{equation*}
N( r,f) =n( 0,f) \cdotp \ln r+\int\limits _{0}^{r}\frac{n( t,f) -n( 0,f)}{t} dt=n( 0,f) \cdotp \ln r+\sum\limits _{0< |t_{k} |< r}\ln\frac{r}{|t_{k} |} \ 
\end{equation*}
 - the average number of poles of $\displaystyle f$ in the disk $\displaystyle |z|\leqslant r$. Now let's denote
\begin{equation*}
m( r,f) =\frac{1}{2\pi }\int\limits _{0}^{2\pi }\ln^{+}\left| f\left( re^{i\varphi }\right)\right| d\varphi \ 
\end{equation*}
- the average value of $\displaystyle \ln |f|\cdot \mathbb{1}_{\{|f| >1\}}$ on the circle $\displaystyle |z|=r$, where $\displaystyle \ln^{+} =\max(1,\ln)$. Let's introduce a Nevanlinna characteristic
\begin{equation*}
T( r,f) =m( r,f) +N( r,f).
\end{equation*}
\end{definition}
The Nevanlinna characteristic for meromorphic functions plays the same role as $\displaystyle M_{r} =\sup _{|z|=r} |f( z) |$ for entire functions: it provides a convenient classification based on growth rate. It's convenient to compare growth rates by comparing them to a power function, measuring its degree in some sense. To define this, let's introduce a more general concept of the order of a function
\begin{gather*}
\rho [ h] =\varlimsup _{r\rightarrow \infty }\frac{\ln^{+} h( r)}{\ln r} \ -\ \text{order} \ h,\\
\lambda [ h] =\varliminf _{r\rightarrow \infty }\frac{\ln^{+} h( r)}{\ln r} \ -\ \text{lower order} \ h.
\end{gather*}
Now let's define the growth orders of a meromorphic function.
\begin{gather*}
\rho \left[ f\right]\xlongequal[\ ]{\mathnormal{def}} \rho \left[ T\left( r,f\right)\right] ,\ \lambda \left[ f\right]\xlongequal[\ ]{\mathnormal{def}} \lambda \left[ T\left( r,f\right)\right],\\
\rho \left[\left\{t_{n}\right\}_{n\in \mathbb{N}}\right]\xlongequal[\ ]{\mathnormal{def}} \rho \left[ n\left( r\right)\right] ,\ \lambda \left[\left\{t_{n}\right\}_{n\in \mathbb{N}}\right]\xlongequal[\ ]{\mathnormal{def}} \lambda \left[ n\left( r\right)\right] ,\ где\ n\left( r\right) =\#\left\{\ k\ |\ |t_{k} |< r\right\}.
\end{gather*}
For the order of sequences, it is known \cite{hayman1964meromorphic} that it is equal to the convergence exponent.
\begin{equation*}
\rho \left[\left\{t_{n}\right\}_{n\in \mathbb{N}}\right] =\inf\left\{\ s\ \Bigl| \ \sum\limits _{n=1}^{\infty }\frac{1}{|t_{n} |^{s}} < +\infty \right\} -\ \text{convergence exponent}.
\end{equation*}
For functions of the type of Cauchy kernel sums $\displaystyle f\left( z\right) =\sum\limits _{k=1}^{\infty }\frac{c_{k}}{z-t_{k}}$, it is known that $\displaystyle \rho \left[ f\right] =\rho \left[\left\{t_{n}\right\}_{n\in \mathbb{N}}\right]$ and $\displaystyle \lambda \left[ f\right] =\lambda \left[\left\{t_{n}\right\}_{n\in \mathbb{N}}\right]$, so their growth is entirely determined by the growth of poles. For functions of finite order, there are results supporting the validity of hypothesis \ref{hyp:ClunieEremenkoRossi} by Clunie, Eremenko, Rossi.

\begin{thm}[M. V. Keldysh, \cite{Keldysh}, 6, p. 159]
Let $\displaystyle \left\{c_{n}\right\}_{n\in \mathbb{N}} \subset \mathbb{C}$ and $\displaystyle \left\{t_{n}\right\}_{n\in \mathbb{N}} \subset \mathbb{C}$, such that $\displaystyle \lim _{n\rightarrow \infty } |t_{n} |=\infty $ and $\displaystyle \left\{t_{n}\right\}_{n\in \mathbb{N}}$ has no finite limit points and
\begin{gather*}
\sum\limits _{n=1}^{\infty } |c_{n} |< +\infty, \\
\sum\limits _{n=1}^{\infty } c_{n} \neq 0,\\
f\left( z\right) =\sum\limits _{k=1}^{\infty }\frac{c_{k}}{z-t_{k}}, \\
\textcolor[rgb]{0.82,0.01,0.11}{\lambda[\left\{t_{n}\right\}_{n\in \mathbb{N}}]<+\infty}.
\end{gather*}
Then $\displaystyle f$ has infinitely many zeros $\displaystyle \left\{s_{n}\right\}_{n\in \mathbb{N}} \subset \mathbb{C}$. Moreover, $\displaystyle \rho \left[\left\{s_{n}\right\}_{n\in \mathbb{N}}\right] = \rho \left[\left\{t_{n}\right\}_{n\in \mathbb{N}}\right]$ and 
\begin{equation*}
1-\delta(0,f)=\varlimsup _{r\rightarrow \infty }\frac{N\left( r,\frac{1}{f}\right)}{T\left( r,f\right)} =1\ -\ \text{"maximum" number of zeros}.
\end{equation*}
\end{thm}
\begin{remark}
In general, the Nevanlinna deficiency of $\displaystyle f$ at the value $\displaystyle a$ is
\begin{equation*}
\delta \left( a,f\right) =\varliminf _{r\rightarrow \infty }\frac{m\left( r,\frac{1}{f-a}\right)}{T\left( r,f\right)}=1-\varlimsup _{r\rightarrow \infty }\frac{N\left( r,\frac{1}{f-a}\right)}{T\left( r,f\right)}.
\end{equation*}
According to the first Nevanlinna theorem (see \cite{hayman1964meromorphic}), $\displaystyle T\left( r,\frac{1}{f-a}\right) =T\left( r,f\right) +O\left( 1\right)$, hence
\begin{equation*}
0\leqslant \delta \left( a,f\right) \leqslant 1.
\end{equation*}
The smaller the deficiency of a function's value, the more "frequently" the function takes it in some sense.
\end{remark}
We will need another decomposition of meromorphic functions, in the form of the quotient of two infinite products. Let's define the Weierstrass factor.
\begin{definition}
\begin{gather*}
E\left( z,p\right) =\left( 1-z\right) e^{z+\frac{z^{2}}{2} +\dotsc +\frac{z^{p}}{p}} ,\ p\in \mathbb{N},\\
E\left( z,0\right) =1-z.
\end{gather*}
\end{definition}
and the genus of the sequence
\begin{definition}
[Genus of a sequence, \cite{Book_Gold_Ostr}]

Let $\displaystyle \left\{t_{n}\right\}_{n\in \mathbb{N}} \subset \mathbb{C}$, then
\begin{equation*}
p\left[\left\{t_{n}\right\}_{n\in \mathbb{N}}\right] =\min\left\{\ p\in \mathbb{N}_{0} \ \Bigl| \ \sum\limits _{n=1}^{\infty }\frac{1}{|t_{n} |^{p+1}} < +\infty \right\} \ -\ \text{genus of the sequence}.
\end{equation*}
\end{definition}
It is known that the genus satisfies the relation
\begin{equation*}
p=\begin{cases}
\left[ \rho \right] , & \rho \notin \mathbb{N}_{0}\\
\left[ \begin{array}{l l}
\rho  & \\
\rho -1 & 
\end{array} \right. & \rho \in \mathbb{N}_{0}
\end{cases}
\end{equation*}
Now let's formulate the Hadamard's theorem
\begin{thm}
[Hadamard, \cite{Book_Gold_Ostr}, \cite{hayman1964meromorphic}]\label{th:hadamard}
Let $\displaystyle f$ be meromorphic in $\displaystyle \mathbb{C}$, $\displaystyle f\left( 0\right) \neq 0$, $\displaystyle \left\{t_{n}\right\}_{n\in \mathbb{N}} ,\left\{s_{n}\right\}_{n\in \mathbb{N}} \subset \mathbb{C}$ - poles and zeros of $\displaystyle f$ with multiplicities and
\begin{gather*}
p\left[\left\{t_{n}\right\}_{n\in \mathbb{N}}\right] \leqslant p,\\
p\left[\left\{s_{n}\right\}_{n\in \mathbb{N}}\right] \leqslant p,\\
\varlimsup _{r\rightarrow \infty }\frac{T\left( r,f\right)}{r^{p+1}} =0.
\end{gather*}
In that case
\begin{equation*}
f\left( z\right) =e^{P\left( z\right)}\frac{\prod\limits _{n} E\left(\frac{z}{s_{n}} ,p\right)}{\prod\limits _{n} E\left(\frac{z}{t_{n}} ,p\right)},
\end{equation*}
where $\displaystyle P\left( z\right)$ is the Taylor polynomial of $\displaystyle \ln f$ of order $\displaystyle p$:
\begin{equation*}
P\left( z\right) =\sum\limits _{k=0}^{p}\left(\ln f\right)^{\left( k\right)}\left( 0\right) \cdotp \frac{z^{k}}{k!}.
\end{equation*}
\end{thm}
It is known that meromorphic functions in $\displaystyle \mathbb{C}$ are the ratios of two entire functions - polynomials of infinite degree, just like rational functions are the ratios of two polynomials. The Hadamard theorem specifies this decomposition, stating that if the meromorphic function does not grow too fast, then it is possible to choose entire functions in the decomposition that do not grow too rapidly either.

\section{The main formula, criterion}

It is known that the Riemann zeta function is defined as
\begin{equation*}
\zeta \left( s\right) =\sum\limits _{n=1}^{\infty }\frac{1}{n^{s}} ,\ \ s >1.
\end{equation*}
It has numerous applications in number theory, particularly in the distribution of prime numbers. Similarly, for an arbitrary sequence $\displaystyle \left\{a_{n}\right\}_{n\in \mathbb{N}}$ of finite order $\displaystyle \rho $, let's define
\begin{equation*}
\zeta _{a}\left( s\right) =\sum\limits _{n=1}^{\infty }\frac{1}{a_{n}^{s}} ,\ \ s >\rho.
\end{equation*}
the zeta function of the sequence
$\displaystyle \left\{a_{n}\right\}_{n\in \mathbb{N}}$. It turns out that for a meromorphic function, the zeta functions of its zeros and poles are related.

\begin{thm}[Connection of zeta functions of zeros and poles of a meromorphic function]\label{th:zetaconnection}

Let $\displaystyle f$ be meromorphic in $\displaystyle \mathbb{C}$, $\displaystyle f(0) \neq 0$, $\displaystyle \lbrace t_{n}\rbrace_{n\in \mathbb{N}} , \lbrace s_{n}\rbrace_{n\in \mathbb{N}} \subset \mathbb{C}$ be the poles and zeros of $\displaystyle f$ with multiplicity, respectively. Let
\begin{gather*}
\sum\limits _{n=1}^{\infty }\frac{1}{|t_{n} |^{N+1}} < +\infty, \\
\sum\limits _{n=1}^{\infty }\frac{1}{|s_{n} |^{N+1}} < +\infty, \\
\varlimsup _{r\rightarrow \infty }\frac{T\left( r,f\right)}{r^{N+1}} =0.
\end{gather*}
Then
\begin{equation*}
\sum\limits _{n=1}^{\infty }\frac{1}{s_{n}^{N+1}} =\sum\limits _{n=1}^{\infty }\frac{1}{t_{n}^{N+1}} -\frac{1}{N!}\left(\frac{f'}{f}\right)^{\left( N\right)}\left( 0\right).
\end{equation*}
\end{thm}
\begin{proof}Given that $\displaystyle \sum\limits _{n=1}^{\infty }\frac{1}{|t_{n} |^{N+1}} < +\infty $ and $\displaystyle \sum\limits _{n=1}^{\infty }\frac{1}{|s_{n} |^{N+1}} < +\infty $, we have $\displaystyle p\left[\left\{t_{n}\right\}_{n\in \mathbb{N}}\right] \leqslant N$ and $\displaystyle p\left[\left\{s_{n}\right\}_{n\in \mathbb{N}}\right] \leqslant N$. Therefore, by Hadamard's theorem \ref{th:hadamard}
\begin{equation*}
f\left( z\right) =e^{P\left( z\right)}\frac{\prod\limits _{n} E\left(\frac{z}{s_{n}} ,N\right)}{\prod\limits _{n} E\left(\frac{z}{t_{n}} ,N\right)},
\end{equation*}
where $\displaystyle P\left( z\right) =\sum\limits _{k=0}^{N}\left(\ln f\right)^{\left( k\right)}\left( 0\right) \cdot \frac{z^{k}}{k!}$ is the Taylor polynomial of $\displaystyle \ln f$ of order $\displaystyle N$. Therefore, for the logarithmic derivative of $\displaystyle f$, we have the relation
\begin{equation*}
\frac{f'\left( z\right)}{f\left( z\right)} =P'\left( z\right) +z^{N}\sum\limits _{n=1}^{\infty }\frac{1}{s_{n}^{N}\left( z-s_{n}\right)} -z^{N}\sum\limits _{n=1}^{\infty }\frac{1}{t_{n}^{N}\left( z-t_{n}\right)}.
\end{equation*}
Let's consider the integral
\begin{equation*}
\frac{1}{2\pi i}\oint\limits _{|z|=r}\frac{1}{z^{N+1}}\frac{f'\left( z\right)}{f\left( z\right)} dz=\frac{1}{N!}\left(\frac{f'}{f}\right)^{\left( N\right)}\left( 0\right) +\sum\limits _{|s_{n} |< r}\frac{1}{s_{n}^{N+1}} -\sum\limits _{|t_{n} |< r}\frac{1}{t_{n}^{N+1}}.
\end{equation*}
Now, let's compute it differently
\begin{gather}
\frac{1}{2\pi i}\oint\limits _{|z|=r}\frac{1}{z^{N+1}}\frac{f'( z)}{f( z)} dz=\frac{1}{2\pi i}\oint\limits _{|z|=r}\frac{P'( z)}{z^{N+1}} dz+\frac{1}{2\pi i}\oint\limits _{|z|=r}\frac{1}{z}\sum\limits _{n=1}^{\infty }\frac{1}{s_{n}^{N}( z-s_{n})} dz-\frac{1}{2\pi i}\oint\limits _{|z|=r}\frac{1}{z}\sum\limits _{n=1}^{\infty }\frac{1}{t_{n}^{N}( z-t_{n})} dz= \notag\\
=\frac{1}{N!}\underbrace{P^{\left( N+1\right)}\left( 0\right)}_{=0} +\frac{1}{2\pi i}\sum\limits _{n=1}^{\infty }\oint\limits _{|z|=r}\frac{1}{z}\frac{1}{s_{n}^{N}\left( z-s_{n}\right)} dz-\frac{1}{2\pi i}\sum\limits _{n=1}^{\infty }\oint\limits _{|z|=r}\frac{1}{z}\frac{1}{t_{n}^{N}\left( z-t_{n}\right)} dz= \notag\\
=\frac{1}{2\pi i}\sum\limits _{n=1}^{\infty }\frac{1}{s_{n}^{N}}\oint\limits _{|z|=r}\frac{1}{z\left( z-s_{n}\right)} dz-\frac{1}{2\pi i}\sum\limits _{n=1}^{\infty }\frac{1}{t_{n}^{N}}\oint\limits _{|z|=r}\frac{1}{z\left( z-t_{n}\right)} dz.
\end{gather}
Let's compute the integrals
\begin{equation*}
\frac{1}{2\pi i}\oint\limits _{|z|=r}\frac{1}{z\left( z-s_{n}\right)} dz=\begin{cases}
\mathop{\mathrm{Res}}\limits_{z = 0}\frac{1}{z\left( z-s_{n}\right)}, & |s_{n} | >r\\
\mathop{\mathrm{Res}}\limits_{z = 0}\frac{1}{z\left( z-s_{n}\right)} +\mathop{\mathrm{Res}}\limits_{z = s_n}\frac{1}{z\left( z-s_{n}\right)}, & |s_{n} |< r
\end{cases} =\begin{cases}
-\frac{1}{s_{n}}, & |s_{n} | >r\\
0, & |s_{n} |< r
\end{cases} =-\frac{\mathbb{1}_{\left\{|s_{n} | >r\right\}}}{s_{n}}.
\end{equation*}
Therefore, we obtain that
\begin{equation*}
\frac{1}{2\pi i}\oint\limits _{|z|=r}\frac{1}{z^{N+1}}\frac{f'\left( z\right)}{f\left( z\right)} dz=-\sum\limits _{|s_{n} | >r}\frac{1}{s_{n}^{N+1}} +\sum\limits _{|t_{n} | >r}\frac{1}{t_{n}^{N+1}}.
\end{equation*}
Comparing, we obtain the statement of the theorem.
\end{proof}
The conditions of the theorem are guaranteed if $\displaystyle N+1 >\rho \left[ f\right]$. Now, let's prove the lemma from the theory of numerical series, which we will need later.
\begin{lema}\label{lem:zero_Series}
Let $\displaystyle \left\{a_{n}\right\}_{n\in \mathbb{N}} \subset \mathbb{C}$ и $\displaystyle \sum\limits\limits _{n=1}^{\infty } |a_{n} |^{M} < +\infty $, then
\begin{equation*}
\left( \ \sum\limits\limits _{n=1}^{\infty } a_{n}^{M} =0,\sum\limits\limits _{n=1}^{\infty } a_{n}^{2M} =0,\sum\limits\limits _{n=1}^{\infty } a_{n}^{3M} =0,\ \dotsc \ \right) \Longrightarrow \left( \forall n\in \mathbb{N} \ \ \ a_{n} =0\right).
\end{equation*}
\end{lema}
\begin{proof}
Then, assuming the contrary, let $\displaystyle \exists n_{0} \in \mathbb{N}$ such that $\displaystyle a_{n_{0}} \neq 0$. Since these series converge, we have $\displaystyle \lim\limits_{n\rightarrow \infty}a_{n}^{M} =0$. Consider the function
\begin{equation*}
f\left( z\right) =\sum\limits _{a_{n} \neq 0}\frac{1}{z-\frac{1}{a_{n}^{M}}}.
\end{equation*}
The series converges uniformly on compacts, allowing for term-by-term differentiation $\displaystyle \infty $ times, yielding a meromorphic function in $\displaystyle \mathbb{C}$. This function exactly has poles at $\displaystyle \frac{1}{a_{n}^{M}}$ for $\displaystyle a_{n} \neq 0$, hence $\displaystyle f\not\equiv 0$. However,
\begin{equation*}
f^{\left( N\right)}\left( 0\right) =-\sum\limits _{a_{n} \neq 0} a_{n}^{M\left( N+1\right)} =0\ \text{при} \ N=0,1,\dotsc .
\end{equation*}
Thus, in the neighborhood of $\displaystyle 0$, $\displaystyle f\left( z\right) =\sum\limits _{N=0}^{\infty } f^{\left( N\right)}\left( 0\right)\frac{z^{N}}{N!} =0$. Consequently, by analytic continuation, we deduce that $\displaystyle f\equiv 0$. This is a contradiction!
\end{proof}

It is interesting to note that venturing into complex analysis and considering the sum of Cauchy kernels provides a concise and elegant proof of this fact, which may initially seem like an 'obvious' lemma from a first-year university course. Now let's examine the formulas relating zeta functions to the zeros and poles of a meromorphic function as stated in Theorem \ref{th:zetaconnection}
\begin{equation*}
\sum\limits _{n=1}^{\infty }\frac{1}{s_{n}^{N+1}} =\sum\limits _{n=1}^{\infty }\frac{1}{t_{n}^{N+1}} -\frac{1}{N!}\left(\frac{f'}{f}\right)^{\left( N\right)}\left( 0\right).
\end{equation*}
If the meromorphic function has no zeros at all, then the relation becomes:
\begin{equation*}
0=\sum\limits _{n=1}^{\infty }\frac{1}{t_{n}^{N+1}} -\frac{1}{N!}\left(\frac{f'}{f}\right)^{\left( N\right)}\left( 0\right) \ \Longrightarrow \sum\limits _{n=1}^{\infty }\frac{1}{t_{n}^{N+1}} =\frac{1}{N!}\left(\frac{f'}{f}\right)^{\left( N\right)}\left( 0\right).
\end{equation*}
It turns out that it is not obvious that this condition is not only necessary for the absence of zeros, but also sufficient.
\begin{thm}
[Criterion for absence of zeros via series]\label{th:krit}
Let $\displaystyle f$ be meromorphic in $\displaystyle \mathbb{C}$, $\displaystyle f\left( 0\right) \neq 0$, $\displaystyle \left\{t_{n}\right\}_{n\in \mathbb{N}} ,\left\{s_{n}\right\}_{n\in \mathbb{N}} \subset \mathbb{C}$ - poles and zeros of $\displaystyle f$ taking multiplicities into account. Let
\begin{gather*}
\sum\limits _{n=1}^{\infty }\frac{1}{|t_{n} |^{N_{0} +1}} < +\infty, \\
\varlimsup _{r\rightarrow \infty }\frac{T\left( r,f\right)}{r^{N_{0} +1}} =0.
\end{gather*}
Then
\begin{equation*}
f\ \text{has no zeros} \ \Longrightarrow \begin{cases}
\sum\limits _{n=1}^{\infty }\frac{1}{t_{n}^{N_{0} +1}} =\frac{1}{N_{0} !}\left(\frac{f'}{f}\right)^{\left( N_{0}\right)}\left( 0\right)\\
\sum\limits _{n=1}^{\infty }\frac{1}{t_{n}^{N_{0} +2}} =\frac{1}{\left( N_{0} +1\right) !}\left(\frac{f'}{f}\right)^{\left( N_{0} +1\right)}\left( 0\right)\\
\dotsc 
\end{cases}
\end{equation*}
Furthermore,
\begin{equation*}
\sum\limits _{n=1}^{\infty }\frac{1}{|s_{n} |^{N_{0} +1}} < +\infty \ \text{и} \ \exists M\geqslant N_{0} +1\ \ \text{т.ч.} \ \ \begin{cases}
\sum\limits _{n=1}^{\infty }\frac{1}{t_{n}^{M}} =\frac{1}{\left( M-1\right) !}\left(\frac{f'}{f}\right)^{\left( M-1\right)}\left( 0\right)\\
\sum\limits _{n=1}^{\infty }\frac{1}{t_{n}^{2M}} =\frac{1}{\left( 2M-1\right) !}\left(\frac{f'}{f}\right)^{\left( 2M-1\right)}\left( 0\right)\\
\dotsc 
\end{cases} \Longrightarrow f\ \text{has no zeros}.
\end{equation*}
\end{thm}
\begin{proof}
1) Let $\displaystyle f$ have no zeros, then $\displaystyle \sum\limits _{n=1}^{\infty }\frac{1}{|s_{n} |^{N_{0} +1}} < +\infty $ is even more true, and then by equations from the theorem \ref{th:krit} we obtain the required.

2) By contradiction. Suppose $\displaystyle f$ has zeros, we know that $\displaystyle \sum\limits _{n=1}^{\infty }\frac{1}{|s_{n} |^{N_{0} +1}} < +\infty $ and $\displaystyle \sum\limits _{n=1}^{\infty }\frac{1}{t_{n}^{kM}} =\frac{1}{N!}\left(\frac{f'}{f}\right)^{\left( kM-1\right)}\left( 0\right)$ for $\displaystyle k\geqslant 1$, then by the equations from theorem \ref{th:krit} we obtain
\begin{equation*}
\sum\limits _{n}\frac{1}{s_{n}^{kM}} =0\ \ \text{at} \ \ k\in \mathbb{N}.
\end{equation*}
Therefore, by Lemma \ref{lem:zero_Series}, we obtain that $\displaystyle \frac{1}{s_{n}^{M}} =0$ for all $\displaystyle n$, whence $\displaystyle s_{n} =\infty $ for all $\displaystyle n$, contradiction!
\end{proof}

In the obtained criterion, the condition $\displaystyle \varlimsup_{r\rightarrow \infty }\frac{T\left( r,f\right)}{r^{N_{0} +1}} =0$ seems difficult, since it is rarely possible to calculate $\displaystyle T\left( r,f\right)$ explicitly, and trivial estimates through the order are often crude and insufficient for applications.

To address this, we present the following classical lemma:
\begin{lema}[\cite{hayman1964meromorphic}]\label{lem:limsup}
Let $\displaystyle \left\{t_{n}\right\}_{n\in \mathbb{N}} \subset \mathbb{C}$ и $\displaystyle n\left( r\right) =\#\left\{\ k\ |\ |t_{k} |< r\right\}$- the counting function, then
\begin{equation*}
\sum\limits _{n=1}^{\infty }\frac{1}{|t_{n} |^{N}} =N\int\limits _{0}^{\infty }\frac{n\left( r\right)}{r^{N+1}} dr=N^{2}\int\limits _{0}^{\infty }\frac{N\left( r\right)}{r^{N+1}} dr.
\end{equation*}
Equality in the sense that if one is finite, then all others are also finite and there are equalities. Additionally,
\begin{equation*}
\sum\limits _{n=1}^{\infty }\frac{1}{|t_{n} |^{N}} < +\infty \ \Longrightarrow \ \varlimsup _{r\rightarrow \infty }\frac{n\left( r\right)}{r^{N}} =\varlimsup _{r\rightarrow \infty }\frac{N\left( r\right)}{r^{N}} =0.
\end{equation*}
\end{lema}
Now, let's simplify the condition in the criterion by removing what was already guaranteed.
\begin{lema}
Let $\displaystyle f$ meromorphic in $\displaystyle \mathbb{C}$, $\displaystyle f\left( 0\right) \neq 0$, $\displaystyle \left\{t_{n}\right\}_{n\in \mathbb{N}} ,\left\{s_{n}\right\}_{n\in \mathbb{N}} \subset \mathbb{C}$ - poles and zeros of $\displaystyle f$ taking multiplicities into account. Let
\begin{gather*}
\sum\limits _{n=1}^{\infty }\frac{1}{|t_{n} |^{N_{0} +1}} < +\infty, \\
\varlimsup _{r\rightarrow \infty }\frac{m\left( r,f\right)}{r^{N_{0} +1}} =0.
\end{gather*}
Then the condition $\displaystyle \varlimsup_{r\rightarrow \infty }\frac{T\left( r,f\right)}{r^{N_{0} +1}} =0$ is satisfied, and thus the conditions of the criterion from Theorem \ref{th:krit} are fulfilled.
\end{lema}
\begin{proof}
Indeed
\begin{equation*}
0\leqslant\varlimsup _{r\rightarrow \infty }\frac{T\left( r,f\right)}{r^{N_{0} +1}} =\varlimsup _{r\rightarrow \infty }\frac{m\left( r,f\right) +N\left( r,f\right)}{r^{N_{0} +1}} \leqslant \underbrace{\varlimsup\limits _{r\rightarrow \infty }\frac{m\left( r,f\right)}{r^{N_{0} +1}}}_{=0\ \text{as stated}} +\underbrace{\varlimsup\limits _{r\rightarrow \infty }\frac{N\left( r,f\right)}{r^{N_{0} +1}}}_{=0\ \text{by lemma \ref{lem:limsup}}} =0
\end{equation*}
\end{proof}
In the future, we will need to apply the criterion to functions of the form of sums of Cauchy kernels and their derivatives, and we need information about the behavior of $\displaystyle m\left( r,f\right)$ to check the condition $\displaystyle \varlimsup\limits _{r\rightarrow \infty }\frac{m\left( r,f\right)}{r^{N_{0} +1}}$. For sums of Cauchy kernels, the Ostrovsky theorem is known:

\begin{thm}
[I. V. Ostrovsky, \cite{Book_Gold_Ostr},6 , p. 252, th. 6.1]
Let $\displaystyle \left\{c_{n}\right\}_{n\in \mathbb{N}} ,\left\{t_{n}\right\}_{n\in \mathbb{N}} \subset \mathbb{C}$, such that $\displaystyle \lim _{n\rightarrow \infty } |t_{n} |=\infty $ и $\displaystyle \left\{t_{n}\right\}_{n\in \mathbb{N}}$ has no finite limit points, and
\begin{gather*}
\sum\limits _{n=1}^{\infty }\frac{|c_{n} |}{|t_{n} |} < +\infty, \\
f\left( z\right) =\sum\limits _{n=1}^{\infty }\frac{c_{n}}{z-t_{n}}.
\end{gather*}
Then for each $\displaystyle 0< p< 1$
\begin{equation*}
\lim _{r\rightarrow \infty }\int\limits _{0}^{2\pi }\left| f\left( re^{i\varphi }\right)\right| ^{p} d\varphi =0.
\end{equation*}
Moreover, $\displaystyle m\left( r,f\right) =o\left( 1\right)$ at $\displaystyle r\rightarrow \infty $.
\end{thm}
For functions of the form of derivatives of sums of Cauchy kernels, I managed to obtain the following theorem.
\begin{thm}
[V. V. Shemyakov, \cite{shemyakov2024infinitesimal}]
Let $\displaystyle \{c_{n}\}_{n\in \mathbb{N}} ,\{t_{n}\}_{n\in \mathbb{N}} \subset \mathbb{C}$, such that $\displaystyle \lim _{n\rightarrow \infty } |t_{n} |=\infty $ and
\begin{equation*}
\{t_{n}\}_{n\in \mathbb{N}} \ \text{has a finite convergence exponent}.
\end{equation*}
Then for each $\displaystyle 0< p< \frac{1}{2}$ and for each $\displaystyle \varepsilon  >0$
\begin{equation*}
\sum\limits _{n=1}^{\infty }\frac{|c_{n} |}{|t_{n} |^{2}} < +\infty \ \Longrightarrow \int\limits _{0}^{2\pi }\left| \sum\limits\limits _{n=1}^{\infty }\frac{c_{n}}{\left( re^{i\varphi } -t_{n}\right)^{2}}\right| ^{p} d\varphi =o\left( r^{\varepsilon }\right) \ \text{at}\ r\rightarrow \infty,
\end{equation*}
\begin{equation*}
\sum\limits _{n=1}^{\infty }\frac{|c_{n} |}{|t_{n} |} < +\infty \ \Longrightarrow \int\limits _{0}^{2\pi }\left| \sum\limits _{n=1}^{\infty }\frac{c_{n}}{\left( re^{i\varphi } -t_{n}\right)^{2}}\right| ^{p} d\varphi =o\left(\frac{1}{r^{p-\varepsilon }}\right) \ \text{at} \ r\rightarrow \infty, 
\end{equation*}
\begin{equation*}
\sum\limits _{n=1}^{\infty } |c_{n} |< +\infty \ \Longrightarrow \int\limits _{0}^{2\pi }\left| \sum\limits _{n=1}^{\infty }\frac{c_{n}}{\left( re^{i\varphi } -t_{n}\right)^{2}}\right| ^{p} d\varphi =o\left(\frac{1}{r^{2p-\varepsilon }}\right) \ \text{at} \ r\rightarrow \infty. 
\end{equation*}
\end{thm}
Now let's formulate the required conditions for the criterion for sums of Cauchy kernels and their derivatives.
\begin{lema}
Let $\displaystyle \left\{c_{n}\right\}_{n\in \mathbb{N}} ,\left\{t_{n}\right\}_{n\in \mathbb{N}} \subset \mathbb{C}$, such that $\displaystyle \lim _{n\rightarrow \infty } |t_{n} |=\infty $ and
\begin{equation*}
f\left( z\right) =\left[ \begin{array}{l l}
\sum\limits _{k=1}^{\infty }\frac{c_{k}}{z-t_{k}} , & \sum\limits _{k=1}^{\infty }\frac{|c_{k} |}{|t_{k} |} < +\infty \\
\sum\limits _{k=1}^{\infty }\frac{c_{k}}{\left( z-t_{k}\right)^{2}} , & \sum\limits _{k=1}^{\infty }\frac{|c_{k} |}{|t_{k} |^{2}} < +\infty 
\end{array} \right.
\end{equation*}
\begin{equation*}
\sum\limits _{n=1}^{\infty }\frac{1}{|t_{n} |^{N_{0} +1}} < +\infty .
\end{equation*}
Then $\displaystyle \varlimsup\limits _{r\rightarrow \infty }\frac{m\left( r,f\right)}{r^{N_{0} +1}} =0$, and the conditions for the absence of zeros from theorem \ref{th:krit} are satisfied.
\end{lema}
\begin{proof}
From the condition $\displaystyle \sum\limits _{n=1}^{\infty }\frac{1}{|t_{n} |^{N_{0} +1}} < +\infty $, we obtain that the order of $\displaystyle \left\{t_{n}\right\}_{n\in \mathbb{N}}$ is finite. Hence, for $\displaystyle \varepsilon =\frac{1}{2}$ in both cases, we have $\displaystyle m\left( r,f\right) =o\left(\sqrt{r}\right)$, from which we have
\begin{equation*}
\varlimsup\limits _{r\rightarrow \infty }\frac{m\left( r,f\right)}{r^{N_{0} +1}} =\varlimsup\limits _{r\rightarrow \infty }\frac{o\left(\sqrt{r}\right)}{r^{N_{0} +1}} = \varlimsup\limits _{r\rightarrow \infty }\frac{o(1)}{r^{N_{0} +\frac{1}{2}}} =0.
\end{equation*}
\end{proof}
Let's examine these systems of equations from the criterion in Theorem \ref{th:krit} in more detail. To do this, we'll define the Bell polynomials.

\begin{definition}
(Bell Polynomials)

Let's define the partial Bell polynomials $\displaystyle \mathcal{B}_{n,k}$ for $\displaystyle 0\leqslant k\leqslant n$.
\begin{equation*}
\mathcal{B}_{n,k}\left( x_{1} ,\dotsc ,x_{n-k+1}\right) =n!\sum\limits _{ \begin{array}{l}
j_{1} +\dotsc +j_{n-k+1} =k\\
1j_{1} +\dotsc +\left( n-k+1\right) j_{n-k+1} =n
\end{array}}\prod _{i=1}^{n-k+1}\frac{x_{i}^{j_{i}}}{\left( i!\right)^{j_{i}} j_{i} !}
\end{equation*}
and complete Bell polynomials
\begin{equation*}
\mathcal{B}_{n}\left( x_{1} ,\dotsc ,x_{n}\right) =\sum\limits _{k=1}^{n}\mathcal{B}_{n,k}\left( x_{1} ,\dotsc ,x_{n-k+1}\right) =n!\sum\limits _{1j_{1} +\dotsc +nj_{n} =n}\prod _{i=1}^{n}\frac{x_{i}^{j_{i}}}{\left( i!\right)^{j_{i}} j_{i} !}.
\end{equation*}
\end{definition}
It is known that Bell polynomials possess an interesting property called "inversion," which is useful for further analysis.

\begin{lema}[Inversion formula for Bell polynomials]\label{lem:inversion}
Let $\displaystyle f\in C^{\infty }\left( a-\varepsilon ;a+\varepsilon \right)$ and $\displaystyle f'\left( a\right) \neq 0$, then
\begin{equation*}
\left( \forall n\in \mathbb{N} \ \ \ y_{n} =\sum\limits _{k=1}^{n} f^{\left( k\right)}\left( a\right)\mathcal{B}_{n,k}\left( x_{1} ,\dotsc ,x_{n-k+1}\right)\right) \Longleftrightarrow \left( \forall n\in \mathbb{N} \ \ \ x_{n} =\sum\limits _{k=1}^{n}\left( f^{-1}\right)^{\left( k\right)}\left( a\right)\mathcal{B}_{n,k}\left( x_{1} ,\dotsc ,x_{n-k+1}\right)\right).
\end{equation*}
\end{lema}

Now let's derive the formula for the $\displaystyle n$-th derivative at zero for a meromorphic function without zeros.

\begin{thm}[Formula for derivatives at zero of a meromorphic function without zeros]
Let $\displaystyle f$ be meromorphic in $\displaystyle \mathbb{C}$, $\displaystyle f\left( 0\right) \neq 0$, $\displaystyle \left\{t_{n}\right\}_{n\in \mathbb{N}} ,\left\{s_{n}\right\}_{n\in \mathbb{N}} \subset \mathbb{C}$ - poles and zeros of $\displaystyle f$ taking multiplicities into account. Let
\begin{gather*}
\sum\limits _{n=1}^{\infty }\frac{1}{|t_{n} |^{N_{0} +1}} < +\infty, \\
\varlimsup _{r\rightarrow \infty }\frac{T\left( r,f\right)}{r^{N_{0} +1}} =0.
\end{gather*}
Then
\begin{equation*}
f\ \text{has no zeros} \ \Longrightarrow \ f^{\left( N\right)}\left( 0\right) =f\left( 0\right) \cdot \mathcal{B}_{N}\left(\left(\ln f\right) '\left( 0\right) ,\dotsc ,\left(\ln f\right)^{\left( N_{0}\right)}\left( 0\right) ,N_{0} !\sum\limits _{n=1}^{\infty }\frac{1}{t_{n}^{N_{0} +1}} ,\dotsc ,\left( N-1\right) !\sum\limits _{n=1}^{\infty }\frac{1}{t_{n}^{N}}\right)
\end{equation*}
\end{thm}
\begin{proof}
Let $\displaystyle f$ have no zeros. Then, according to the criterion from Theorem \ref{th:krit} of zero absence, we obtain an infinite system of equations
\begin{equation*}
T_{N} =\sum\limits _{n=1}^{\infty }\frac{1}{t_{n}^{N}} =\frac{1}{\left( N-1\right) !}\left(\frac{f'}{f}\right)^{\left( N-1\right)}\left( 0\right) \ \text{при} \ N\geqslant N_{0} +1.
\end{equation*}
Let's extend $\displaystyle T_{k} =\frac{1}{\left( k-1\right) !}\left(\frac{f'}{f}\right)^{\left( k-1\right)}\left( 0\right)$ for $\displaystyle k=1..N_{0}$ for convenience, then
\begin{equation*}
T_{N} =\frac{1}{\left( N-1\right) !}\left(\frac{f'}{f}\right)^{\left( N-1\right)}\left( 0\right) \ \ \text{при} \ N\geqslant 1.
\end{equation*}
Now, let's apply Faà di Bruno's formula and expand the multiple derivative using Bell polynomials
\begin{equation*}
\left( N-1\right) !T_{N} =\left(\frac{f'}{f}\right)^{\left( N-1\right)}\left( 0\right) =\left(\ln f\right)^{\left( N\right)}\left( 0\right) =\sum\limits _{k=1}^{N}\ln^{\left( k\right)}\left( f\left( 0\right)\right) \cdot \mathcal{B}_{N,k}\left( f'\left( 0\right) ,\dotsc ,f^{\left( N+1-k\right)}\left( 0\right)\right).
\end{equation*}
Therefore, according to the inversion formula for Bell polynomials from Lemma \ref{lem:inversion}, we obtain
\begin{gather*}
f^{\left( N\right)}\left( 0\right) =\sum\limits _{k=1}^{N}\left( e^{z}\right)^{\left( k\right)}\left(\ln f\left( 0\right)\right) \cdot \mathcal{B}_{N,k}\left( T_{1} ,\dotsc ,\left( N-k\right) !T_{N+1-k}\right) =\\
=f\left( 0\right)\sum\limits _{k=1}^{N}\mathcal{B}_{N,k}\left( T_{1} ,\dotsc ,\left( N-k\right) !T_{N+1-k}\right) =f\left( 0\right) \cdot \mathcal{B}_{N}\left( T_{1} ,\dotsc ,\left( N-1\right) !T_{N}\right).
\end{gather*}
\end{proof}
From this formula, it follows that for a meromorphic function $\displaystyle f$ of finite order without zeros, knowing all the poles and a sufficient number of derivatives of $\displaystyle \ln f$ at zero completely determines the function $\displaystyle f$.
 
\section{Application of the Criterion}

Let us consider several applications of the zero absence criterion.

\begin{thm}
[Existence of a zero for sums of Cauchy kernels with residues 1]
Let $\displaystyle \{t_{n}\}_{n\in \mathbb{N}} \subset \mathbb{C}$, such that $\displaystyle \lim _{n\rightarrow \infty } |t_{n} |=\infty $ и , and
\begin{gather*}
\sum _{n=1}^{\infty }\frac{1}{|t_{n} |} < +\infty, \\
\sum _{n=1}^{\infty }\frac{1}{t_{n}} \neq 0.
\end{gather*}
Then the meromorphic function $\displaystyle f( z) =\sum _{k=1}^{\infty }\frac{1}{z-t_{k}}$ \textcolor[rgb]{0.82,0.01,0.11}{has zeros} at $\displaystyle \mathbb{C}$.
\end{thm}
\begin{proof}
By contradiction, suppose $\displaystyle f$ has no zeros. Then, according to the criterion from Theorem\ref{th:krit}
\begin{equation*}
\sum\limits _{n=1}^{\infty }\frac{1}{t_{n}^{N+1}} =\frac{1}{N!}\left(\frac{f'}{f}\right)^{\left( N\right)}\left( 0\right) \ \text{при} \ N\geqslant N_{0} =0.
\end{equation*}
Notice that $\displaystyle \left(f\right)^{(N)}(0) = -N!\sum_{n=1}^{\infty} \frac{1}{t_n^{N+1}}$. Therefore, we obtain the relation
\begin{equation*}
\left(\frac{f'}{f} +f\right)^{\left( N\right)}\left( 0\right) =0\ \text{при} \ N\geqslant 0.
\end{equation*}
Then, by the uniqueness theorem for analytic functions, we obtain a differential equation in some neighborhood of $\displaystyle 0$
\begin{equation*}
\frac{f'}{f} +f=0\ \Longleftrightarrow \ f( z) =\frac{1}{z-C}.
\end{equation*}
Но $\displaystyle f( z)$ cannot be rational, a contradiction!
\end{proof}

This theorem actually follows from the results of Clunie, Eremenko, Rossi (\cite{ClunieEremenkoRossi}, th. 2.1),(\cite{LangleyRossiZeros}, th. 1.3) but the proof method itself is interesting. Now let's obtain the theorem for derivatives of sums of Cauchy kernels.

\begin{thm}
[Derivatives of sums of Cauchy kernels with residues 1 without zeros]
Let $\displaystyle \{t_{n}\}_{n\in \mathbb{N}} \subset \mathbb{C}$, such that $\displaystyle \lim _{n\rightarrow \infty } |t_{n} |=\infty $ и , and
\begin{gather*}
\sum _{n=1}^{\infty }\frac{1}{|t_{n} |^{2}} < +\infty \ -\ \text{natural condition},\\
\sum _{n=1}^{\infty }\frac{1}{t_{n}^{2}} \neq 0\ -\ \text{restriction}.
\end{gather*}
Then
\begin{equation*}
f( z) =\sum _{k=1}^{\infty }\frac{1}{( z-t_{k})^{2}} \ \textcolor[rgb]{0.82,0.01,0.11}{\text{has no zeros}} \ \Longrightarrow \ f( z) =\frac{\pi ^{2}}{\sin^{2}( \pi z+C)} ,\ C\in \mathbb{C}.
\end{equation*}
\end{thm}
\begin{proof}
By contradiction, suppose $\displaystyle f$ has no zeros. Then, according to the criterion from theorem \ref{th:krit}
\begin{equation*}
\sum\limits _{n=1}^{\infty }\frac{2}{t_{n}^{N+1}} =\frac{1}{N!}\left(\frac{f'}{f}\right)^{( N)}( 0) \ \text{at} \ N\geqslant N_{0} =1.
\end{equation*}
Notice that $\displaystyle (f)^{(N-1)}(0) = N!\sum_{n=1}^{\infty} \frac{1}{t_n^{N+1}}$. Therefore, we obtain the relation
\begin{equation*}
\left(\left(\frac{f'}{f}\right)^{'} -2f\right)^{( N-1)}( 0) =0\ \text{при} \ N\geqslant 1.
\end{equation*}
Therefore, by the uniqueness theorem for analytic functions, we obtain a differential equation in a neighborhood of $\displaystyle 0$
\begin{equation*}
\left(\frac{f'}{f}\right)^{'} -2f=0\ \Longleftrightarrow \ f( z) =\frac{1}{z-C} \ \text{или} \ f( z) =\frac{\pi ^{2}}{\sin^{2}( \pi z+C)}.
\end{equation*}
But $\displaystyle f(z)$ cannot be rational, so $\displaystyle f(z) = \frac{\pi^2}{\sin^2(\pi z + C)}$.
\end{proof}

In fact, we can significantly weaken the requirements in the theorem while involving a much more powerful theory, for example see (\cite{baranov2023zeros}, p. 11, th. 1.4).

Now let's give geometric conditions on poles for the existence of zeros of sums of Cauchy kernels. We will require, as in the Clunie, Eremenko, Rossi hypothesis \ref{hyp:ClunieEremenkoRossi}, the condition $\displaystyle c_{n} > 0$. The following two theorems \ref{thm:pi/4} and \ref{thm:pi/8} (a consequence of theorem \ref{th:krit}) were obtained by my science supervisor Anton Baranov.
\begin{thm}\label{thm:pi/4}
Let $\displaystyle \{t_{n}\}_{n\in \mathbb{N}} \subset \mathbb{C}$, such that $\displaystyle \lim _{n\rightarrow \infty } |t_{n} |=\infty $ and
\begin{gather*}
c_{n}> 0,\ | \arg \ t_{n}| < \frac{\pi }{4} \ \ \forall n\in \mathbb{N}, \\
\sum _{n=1}^{\infty }\frac{c_{n}}{|t_{n} |} < +\infty ,\\
\sum _{n=1}^{\infty }\frac{1}{|t_{n} |} < +\infty.
\end{gather*}
Then the meromorphic function $\displaystyle f(z) = \sum_{n=1}^{\infty} \frac{c_{n}}{z-t_{n}}$  \textcolor[rgb]{0.82,0.01,0.11}{\text{has zeros}} at $\mathbb{C}$.
\end{thm}
\begin{proof}
By contradiction, suppose $\displaystyle f$ has no zeros. Then, according to the criterion from theorem \ref{th:krit}
\begin{equation*}
\sum\limits _{n=1}^{\infty }\frac{1}{t_{n}^{N+1}} =\frac{1}{N!}\left(\frac{f'}{f}\right)^{( N)}( 0) \ \text{при} \ N\geqslant N_{0} =0.
\end{equation*}
Therefore, we obtain the relation
\begin{equation*}
\sum\limits _{n=1}^{\infty }\frac{1}{t_{n}} =\frac{f'( 0)}{f( 0)} =\frac{\sum\limits\limits _{n=1}^{\infty }\frac{c_{n}}{t_{n}^{2}}}{\sum\limits\limits _{n=1}^{\infty }\frac{c_{n}}{t_{n}}}.
\end{equation*}
Then
\begin{equation*}
\sum\limits\limits _{n=1}^{\infty }\frac{c_{n}}{t_{n}^{2}} =\sum\limits _{n=1}^{\infty }\frac{1}{t_{n}} \cdot \sum\limits\limits _{n=1}^{\infty }\frac{c_{n}}{t_{n}} =\sum\limits\limits _{n=1}^{\infty }\frac{c_{n}}{t_{n}^{2}} +\sum\limits _{n=1}^{\infty }\frac{1}{t_{n}}\sum\limits _{k\neq n}\frac{c_{k}}{t_{n}}.
\end{equation*}
From here, we get that
\begin{equation*}
\sum\limits _{n=1}^{\infty }\frac{1}{t_{n}}\sum\limits _{k\neq n}\frac{c_{k}}{t_{n}} =0.
\end{equation*}
But this cannot be, as
\begin{equation*}
0=\Re \sum\limits _{n=1}^{\infty }\frac{1}{t_{n}}\sum\limits _{k\neq n}\frac{c_{k}}{t_{n}} =\sum\limits _{n=1}^{\infty }\sum\limits _{k\neq n}\underbrace{c_{k}}_{ >0}\underbrace{\Re \left(\frac{1}{t_{k} t_{n}}\right)}_{ >0}  >0.
\end{equation*}
Contradiction!
\end{proof}

\begin{thm}\label{thm:pi/8}
Let $\displaystyle \{t_{n}\}_{n\in \mathbb{N}} \subset \mathbb{C}$, such that $\displaystyle \lim _{n\rightarrow \infty } |t_{n} |=\infty $ and
\begin{gather*}
c_{n}> 0,\ | \arg \ t_{n}| < \frac{\pi }{8} \ \ \forall n\in \mathbb{N}, \\
\sum _{n=1}^{\infty }\frac{c_{n}}{|t_{n}|} < +\infty ,\\
\sum _{n=1}^{\infty }\frac{1}{|t_{n}|^2} < +\infty.
\end{gather*}
Then the meromorphic function $\displaystyle f(z) = \sum_{n=1}^{\infty} \frac{c_{n}}{z-t_{n}}$  \textcolor[rgb]{0.82,0.01,0.11}{\text{has zeros}} at $\mathbb{C}$.
\end{thm}
\begin{proof}
 By contradiction, suppose $\displaystyle f$ has no zeros. Then, according to the criterion from theorem \ref{th:krit}
\begin{equation*}
\sum\limits _{n=1}^{\infty }\frac{1}{t_{n}^{N+1}} =\frac{1}{N!}\left(\frac{f'}{f}\right)^{( N)}( 0) \ \text{при} \ N\geqslant N_{0} =1.
\end{equation*}
Therefore, we obtain the relation
\begin{equation*}
\sum\limits _{n=1}^{\infty }\frac{1}{t_{n}^{2}} =\left(\frac{f'}{f}\right)^{'}( 0) =\frac{f''( 0) f( 0) -f'( 0)^{2}}{f^{2}( 0)}
\end{equation*}.
Then
\begin{equation*}
\sum\limits\limits _{n=1}^{\infty }\frac{1}{t_{n}^{2}} \cdot \left(\sum\limits\limits _{n=1}^{\infty }\frac{c_{n}}{t_{n}}\right)^{2} =2\sum\limits\limits _{n=1}^{\infty }\frac{c_{n}}{t_{n}} \cdot \sum\limits\limits _{n=1}^{\infty }\frac{c_{n}}{t_{n}^{3}} -\left(\sum\limits\limits _{n=1}^{\infty }\frac{c_{n}}{t_{n}^{2}}\right)^{2}.
\end{equation*}
Let $\displaystyle \beta =\sum_{n=1}^{\infty }\frac{c_{n}}{t_{n}}$. Then
\begin{gather*}
0=\sum\limits\limits _{n=1}^{\infty }\left(\frac{\beta ^{2}}{t_{n}^{2}} -2\frac{\beta c_{n}}{t_{n}^{3}} +\frac{c_{n}^{2}}{t_{n}^{4}} +\frac{c_{n}}{t_{n}^{2}}\sum\limits\limits _{k\neq n}\frac{c_{k}}{t_{k}^{2}}\right) =\sum\limits\limits _{n=1}^{\infty }\frac{1}{t_{n}^{2}}\left(\left( \beta -\frac{c_{n}}{t_{n}}\right)^{2} +c_{n}\sum\limits\limits _{k\neq n}\frac{c_{k}}{t_{k}^{2}}\right) =\\
=\underbrace{\sum\limits\limits _{n=1}^{\infty }\underbrace{\frac{1}{t_{n}^{2}}}_{|\arg |< \pi /4}\left(\underbrace{\left(\sum\limits\limits _{k\neq n}\frac{c_{k}}{t_{k}}\right)^{2}}_{|\arg |< \pi /4} +\underbrace{c_{n}\sum\limits\limits _{k\neq n}\frac{c_{k}}{t_{k}^{2}}}_{|\arg |< \pi /4}\right)}_{\Re  >0}.
\end{gather*}
Contradiction!
\end{proof}


\begin{thebibliography}{CER93}

\bibitem[AE94]{LangleyRossiZeros}
J.~Rossi A.~Eremenko, J.~Langley.
\newblock On the zeros of meromorphic functions of the form $\sum\limits_n \dfrac{c_n}{z-t_n}$.
\newblock {\em Journal d’Analyse Mathematique}, 1994.

\bibitem[AG08]{Book_Gold_Ostr}
I.V.~Ostrovskii A.A.~Goldberg.
\newblock {\em Value distribution of meromorphic functions}, volume 236.
\newblock Translations of Mathematical Monographs, 2008.

\bibitem[BS23]{baranov2023zeros}
Anton Baranov and Vladimir Shemyakov.
\newblock Zeros of meromorphic functions of the form $\sum\limits_n \dfrac{c_n}{(z-t_n)^2}$.
\newblock {\em conditinally accepted to J. Math. Anal. Appl.}, 2023.
\newblock arXiv:2307.04211 [math.CV].

\bibitem[CER93]{ClunieEremenkoRossi}
J.~Clunie, A.~Eremenko, and J.~Rossi.
\newblock On equilibrium points of logarithmic and newtonian potentials.
\newblock {\em Journal of the London Mathematical Society}, s2-47(2):309--320, 1993.

\bibitem[Hay64]{hayman1964meromorphic}
W.K. Hayman.
\newblock {\em Meromorphic Functions}.
\newblock Oxford mathematical monographs. Clarendon Press, 1964.

\bibitem[Kel85]{Keldysh}
M.~V. Keldysh.
\newblock {\em {Izbrannye trudy. Matematika}}.
\newblock ``Nauka'', Moscow, 1985.
\newblock Edited and with a preface by K. I. Babenko, N. N. Bogolyubov, and N. N. Chentsov.

\bibitem[She24]{shemyakov2024infinitesimal}
Vladimir Shemyakov.
\newblock The infinitesimal behavior of the sum of cauchy kernels and its derivative at infinity.
\newblock 2024.
\newblock arXiv:2403.04152 [math.CV].

\end{thebibliography}
\end{document}